\documentclass[10pt,reqno]{amsart}
\usepackage{amssymb, amsmath, euscript, verbatim, default}
\usepackage{anisotropic}
\newcommand{\N}{\mathbb{N}}
\newcommand{\Z}{\mathbb{Z}}
\newcommand{\R}{\mathbb{R}}

\renewcommand{\theenumi}{(\roman{enumi})}

\newcommand{\eps}{\epsilon}
\newcommand{\X}{\mathcal{X}}
\newcommand{\Y}{\mathcal{Y}}

\begin{document}

\title[Anisotropic Function Spaces]{Anisotropic Function Spaces and Elliptic\\Boundary Value Problems}
\author{Timothy Nguyen}
\date{\today}

\begin{abstract}
  In this paper, we study anisotropic Bessel potential and Besov spaces, where the anisotropy measures the extra amount of regularity in certain directions.  Some basic properties of these spaces are established along with applications to elliptic boundary value problems.
\end{abstract}

\maketitle

\section{Introduction}

The theory of elliptic boundary value problems is well understood on a variety of function spaces.  Among these function spaces, the ones of particular importance include the classical Besov spaces and Triebel-Lizorkin spaces (see \cite{Tr}, \cite{RS})

In this paper, we collect some basic properties of certain anisotropic function spaces.  The modifier anisotropic means that the function space captures different amounts of regularity in different directions.  However, there is great freedom in how one may choose to define the anisotropy.  Frequently in the literature on function spaces, the anisotropy measures how different directions scale with respect to one another (e.g., in parabolic problems, there is a two to one ratio between space and time).  However for us, the anisotropy will measure how much additional regularity one has in certain directions (see Definition \ref{Def}), and thus is of an additive nature instead of a multiplicative one.  Such spaces occur, for instance, in \cite{Ho} in the context of elliptic boundary value problems, where $L^2$ type spaces are used.  Our paper will study the Besov and Bessel potential generalization of these anisotropic spaces.

In studying the general basic properties of these anisotropic spaces, most of them are direct consequences of their classical (isotropic) counterparts, since in many instances, the anisotropy simply carries through with little or no modification.  In some cases however, one has to work harder and for this, one can apply results from product type spaces, also called spaces of dominating mixed smoothness (see \cite{ST}).  These latter spaces are also anisotropic spaces in the broad sense as described above and are closely related to the anisotropic spaces modeled on those of H\"ormander which we are interested in studying.  By results of Yamazaki in \cite{Ya1} and \cite{Ya2}, we will see that pseudodifferential operators and product type pseudodifferential operators are bounded both on spaces of dominating mixed smoothness and on our anisotropic function spaces.  This fact will help us greatly in proving some of the basic properties about our anisotropic function spaces.  Finally, we apply all these results to generalize the standard elliptic estimates for elliptic boundary value problems on the classical (isotropic) function spaces to our anisotropic function spaces.

The author's motivation for this work arises from the occurrence of anisotropic function spaces in the study of the nonlinear PDE in \cite{N2}. Thus, the results compiled here comprise a balance between the author's needs and sufficient generality.  

\section{Definitions}

We begin by defining our function spaces on Euclidean space.  From this, it is routine to define these spaces for subsets of Euclidean space and for manifolds (with boundary).

On $\R^n$, let $\S$ denote the space of rapidly decaying Schwartz functions and let $\S'$ denote its dual space of tempered distributions.  Let $\F$ denote the Fourier transform
\begin{equation}
 \F f(\xi) = \int e^{i\xi\cdot x}f(x)dx
\end{equation}
acting on tempered distributions $\S'$.  We will also write $\hat f = \F f$ for shorthand. Let $\{\varphi_j\}_{j=0}^\infty$ be a dyadic Littlewood-Paley partition of unity of $\R^n$.  In other words, let $\varphi_0$ be a smooth compactly supported bump function with $\varphi_0 \leq 1$, $\varphi_0(x) = 1$ for $|x| \leq 1$, and $\varphi_0(x) = 0$ for $|x| \geq 2$.  From this, define $\varphi_1(x) = \varphi_0(x/2) - \varphi_0(x)$ and $\varphi_j(x) = \varphi_1(2^{1-j}x)$ for all $j \geq 2$.  From this, we obtain a dyadic partition of unity, since $\supp \varphi_j \subseteq [2^{j-1},2^{j+1}]$ for $j \geq 1$ and $\sum_{j=0}^\infty \varphi_j \equiv 1$.

Given $f \in \S'$, let
$$f_j = \F^{-1}\varphi_j\F f$$
be the $j$th dyadic piece of $f$.  Then $f = \sum_{j=0}^\infty f_j$, and from this Littlewood-Paley decomposition, we recall the following definition of the classical Triebel-Lizorkin and Besov spaces:

\begin{Definition}Let $s \in \R$ and $0<p,q<\infty$.
 \begin{enumerate}
  \item Define the \textit{Triebel-Lizorkin spaces}
\begin{equation}
 F^s_{p,q}(\R^n) := \{f \in \S'(\R^n): \|f\|_{F^s_{p,q}}\| = \|\{2^{sj}f_j\}_{\ell^q_j}\|_{L^p} < \infty\}.
\end{equation}
\item Define the \textit{Besov spaces}
\begin{equation}
 B^s_{p,q}(\R^n) := \{f \in \S'(\R^n) : \|f\|_{B^s_{p,q}} = \{\|2^{sj}f_j\|_{L^p}\}_{\ell^q_j} < \infty\}.
\end{equation}
 \end{enumerate}
\end{Definition}

In this paper, for simplicity, we will restrict ourselves to the most frequently occurring of these spaces, at least, in the context of elliptic boundary value problems.  Namely, we restrict ourselves to the range $1<p<\infty$, and we will work with the spaces
\begin{align}
 H^{s,p}(\R^n) := F^s_{p,2}(\R^n)\\
 B^{s,p}(\R^n) := B^s_{p,p}(\R^n).
\end{align}
As is well known, when $s$ is a nonnegative integer, $H^{s,p}(\R^n)$ is the usual Sobolev space of functions whose derivatives up to order $s$ lie in $L^p$.  Moreover, for $s > 1/p$, the space $B^{s-1/p,p}(\R^{n-1})$ is the image of $H^{s,p}(\R^n)$ under the restriction map from $\R^n$ to $\R^{n-1}$.  Furthermore, we have $H^{s,2}=B^{s,2}$ for all $s \in \R$.  The spaces $H^{s,p}(\R^n)$ are also called \textit{Bessel potential spaces}, and we will refer to them as such.

We now adapt the above construction to define product type spaces.  Suppose we have a product decomposition
$$\R^n = \R^{n_1} \times \R^{n_2}.$$
Let $x^{(1)}$ and $x^{(2)}$ be the variables on $\R^{n_1}$ and $\R^{n_2}$, respectively, and likewise, let $\xi^{(1)}$ and $\xi^{(2)}$ be the corresponding dual Fourier variables.  As before, let $\{\varphi_j\}_{j=0}^\infty$ denote the dyadic partition of unity on $\R^{n_1}\times\R^{n_2}$ defined as above, and similarly, let $\{\varphi^{(i)}_j\}_{j=0}^\infty$ be the corresponding dyadic partition of unity on $\R^{n_i}$, $i = 1,2$.  Thus, in addition to the usual radial Littlewood-Paley decomposition on $\R^{n_1}\times\R^{n_2}$, we also get a product Littlewood-Paley decomposition for $f \in \S'(\R^{n_1}\times\R^{n_2})$, namely
$$f = \sum_{j,k \geq 0} f_{j,k}$$
where
$$f_{j,k} = \F^{-1}\varphi^{(1)}_j\varphi^{(2)}_k\F f.$$

From this product decomposition, one can define product-type spaces.

\begin{Definition}\label{DefProdSpaces} Let $s,t\in\R$ and $0 <p,q <\infty$.  Then define the \textit{product type Triebel Lizorkin spaces} and \textit{product type Besov spaces}
\begin{align}
 F^{s,t}_{p,q}(\R^{n_1}\times\R^{n_2}) &= \{f \in \S'(\R^n) : \|f\|_{F^{s,t}_{p,q}} = \|\{2^{js+kt}f_{j,k}\}_{\ell^q_{j,k}}\|_{L^p}\}\\
 B^{s,t}_{p,q}(\R^{n_1}\times\R^{n_2}) &= \{f \in \S'(\R^n) : \|f\|_{F^{s,t}_{p,q}} = \{\|2^{js+kt}f_{j,k}\|_{L^p}\}_{\ell^q_{j,k}}\},
\end{align}
respectively.
\end{Definition}

These spaces are also known as spaces of dominating mixed smoothness (see e.g. \cite{ST}).They are also anisotropic in the sense that different directions have different amounts of regularity.  However, these particular spaces will only play an auxiliary role in what we do.  The type of anisotropic function spaces we will be defining are those that are obtained from the classical Bessel potential and Besov spaces by specifying an extra degree of smoothness in the second factor $\R^{n_2}$.  To measure this, we introduce the following family of Bessel potential operators acting on $\R^{n_2}$
\begin{equation}
 J_{(2)}^{s_2}f = \F^{-1}\big<\xi^{(2)}\big>^{s_2}\F f, \qquad s_2 \in \R,
\end{equation}
where $\left<\xi^{(2)}\right> := (1 + |\xi^{(2)}|^2)^{1/2}$ for $\xi^{(2)} \in \R^{n_2}$.  These operators constitute the anisotropic version of the usual Bessel potential operators acting on $\R^n$:
\begin{equation}
 J^sf = \F^{-1}\left<\xi\right>^s\F f, \qquad s \in \R.
\end{equation}

\begin{Definition}\label{Def} Let $s_1,s_2 \in \R$, and $1<p<\infty$.
 \begin{enumerate}
  \item Define the \textit{anisotropic Bessel potential spaces}
\begin{equation}
 H^{(s_1,s_2),p}(\R^{n_1}\times\R^{n_2}) = \{f \in S'(\R^n): \|f\|_{H^{(s_1,s_2),p}} = \|J_{(2)}^{s_2}f\|_{H^{s_1,p}} < \infty\}.
\end{equation}
\item Define the \textit{anisotropic Besov spaces}
\begin{equation}
 B^{(s_1,s_2),p}(\R^{n_1}\times\R^{n_2}) = \{f \in S'(\R^n): \|f\|_{B^{(s_1,s_2),p}} = \|J_{(2)}^{s_2}f\|_{B^{s_1,p}} < \infty\}.
\end{equation}
 \end{enumerate}
\end{Definition}

\noindent\textbf{Notation. }As an abbreviation, we write $A^{(s_1,s_2),p}$ as shorthand for either $H^{(s_1,s_2),p}$ or $B^{(s_1,s_2),p}$.  If a formula appears with multiple occurences of $A^{(\bullet,\bullet),\bullet}$, then we always mean that all the spaces are simultaneously anisotropic Bessel potential or anisotropic Besov spaces.  Similarly, we also write $A^{s,p} = A^{(s,0),p}$ as shorthand for $H^{s,p}$ or $B^{s,p}$.\\

Thus the $s_2$ parameter in the above definitions measures the extra degree of differentiability in the $\R^{n_2}$ direction.  From now on, when we refer to anisotropic function spaces, we will mean the spaces $A^{(s_1,s_2),p}$.  When $s_2=0$ we will say the space is isotropic.
We have the following basic relations among our anisotropic spaces:

\begin{Lemma}\label{LemmaLift}
Let $s_1,s_2 \in \R$, and $1<p<\infty$.
\begin{enumerate}
  \item (Lift Property) The map $J^{s_1'}J_{(2)}^{s_2'}: A^{(s_1,s_2),p}(\R^{n_1} \times \R^{n_2})  \to A^{(s_1-s_1',s_2-s_2'),p}(\R^{n_1} \times \R^{n_2})$ is an isomorphism for all $s_1',s_2' \in \R$.
   \item Let $s_2 \geq 0$ and let $D_{s_2}$ be any elliptic pseudodifferential operator of order $s_2$ on $\R^{n_2}$.  Then we have the following equivalence of norms
  \begin{equation}
    \frac{1}{C}\|f\|_{A^{(s_1,s_2),p}} \leq \|f\|_{A^{s_1,p}} + \|D_{s_2} f\|_{A^{s_1,p}} \leq C\|f\|_{A^{(s_1,s_2),p}},
  \end{equation}
  where $C$ depends only on $s_1$,$s_2$,and $p$.  
\end{enumerate}
\end{Lemma}

\begin{proof} (i) This follows from the definition of $A^{(s_1,s_2),p}$ and the usual lift property of the isotropic spaces $A^{s,p}$, which states that $J^t: A^{s,p}(\R^n) \to A^{s-t,p}(\R^n)$ is an isomorphism for all $s,t \in \R$ (see \cite{Tr}). (ii) This is a simple consequence of the fact that product type pseudodifferential operators are bounded on the spaces $A^{(s_1,s_2),p}$, which we prove in the next section.
\end{proof}

We will provide some other equivalent descriptions of the spaces $A^{(s_1,s_2),p}$ in Lemma \ref{LemmaEquiv}. In the remaining sections of this paper, we work out some of the basic properties of these anisotropic function spaces, which are modeled on the basic properties that their isotropic counterparts satisfy.  We then apply these results to the study of elliptic boundary value problems in the final section.

For now, we record the following important fact:

\begin{Lemma}\label{LemmaProdLP}
  For all $1<p<\infty$, we have $F^{0,0}_{p,2}(\R^{n_1}\times\R^{n_2}) = F^0_{p,2}(\R^n)$.
\end{Lemma}

In other words, the lemma tells us that in addition to the classical radial Littlewood-Paley theorem, which states $F^0_{p,2} = L^p$, we also have a product type Littlewood-Paley theorem, since the lemma implies $F^{0,0}_{p,2} = L^p$.  This fact will be important in the proof of Theorem \ref{ThmProdPSDO}. A proof of Lemma \ref{LemmaProdLP} can be found in e.g. \cite{Ya1}.

\section{Mapping Properties of Pseudodifferential Operators}

It is well-known that pseudodifferential operators are bounded on Triebel-Lizorkin and Besov spaces.  What is less well-known is that product type pseudodifferential operators are also bounded on these spaces and their product type counterparts.  For product type pseudodifferential operators that are purely Fourier multipliers, the result can be found in \cite{ST}, and the general case can be found in \cite{Ya2}. As we will show, it readily follows from these results that product type pseudodifferential operators are bounded on anisotropic spaces.

Let us recall the definition of a pseudodifferential operator on $\R^n$ so that we may define precisely what a product type pseudodifferential operator is.  For every $m \in \R$, we can define the symbol class $S^m = S^m(\R^n)$ to be the space of all smooth functions $a(x,\xi) \in C^\infty(\R^n \times \R^n)$ such that
\begin{equation}
 \sup_{x,\xi}|\partial_x^\beta\partial_\xi^\alpha a(x,\xi)| \leq C_{\alpha,\beta}\left<\xi\right>^{m-|\alpha|}
\end{equation}
for all multi-indices $\alpha$,$\beta$.  Thus, the space $S^m$ is a Fr\'echet space whose topology is generated by the semi-norms
\begin{equation}
\|a\|_{S^m_{\alpha,\beta}} := \sup_{x,\xi}\left<\xi\right>^{-m+|\alpha|}|a(x,\xi)|.
\end{equation}
Given a symbol $a(x,\xi) \in S^m$, we obtain the associated $m$th order pseudodifferential operator
\begin{equation}
 a(x,D)f = (2\pi)^{-n}\int e^{i(x-y)\cdot\xi}a(x,\xi)f(y)dyd\xi, \label{quantize}
\end{equation}
defined for $f \in C^\infty_0(\R^n)$.  Here, $D = (i^{-1}\partial_{x^1}, \ldots, i^{-1}\partial_{x^n})$ is the operator given by the Fourier multiplier $\xi$.  Let $OS^m$ denote the class of all $m$th order pseudodifferential operators obtained by (\ref{quantize}) for $a \in S^m$.

We have the following standard theorem concerning pseudodifferential operators:

\begin{Theorem}\label{ThmPSDO}\hspace{.1in}
 \begin{enumerate}
  \item For all $m_1,m_2 \in \R$, we have the composition rule $OS^{m_1} \circ OS^{m_2} \to OS^{m_1+m_2}$.
  \item If $P \in OS^0$, then $P$ is bounded on $A^{s,p}(\R^n)$ for all $1<p<\infty$ and $s \in \R$.  Moreover, for any fixed $s$ and $p$, the operator norm of $P$ is bounded in terms of only finitely many symbol semi-norms $S^0_{\alpha,\beta}$.
 \end{enumerate}
\end{Theorem}

Next, we will define some more general symbol classes in order to generalize the above theorem to anisotropic function spaces and to a wider class of operators.  Suppose we have a decomposition $\R^n = \R^{n_1}\times\R^{n_2}$.  As before, write $x,\xi \in \R^n$ as $(x^{(1)},x^{(2)})$ and $(\xi^{(1)},\xi^{(2)})$ with respect to this decomposition.  Likewise, if $\alpha \in \Z_+^n$ is a multi-index of nonnegative integers, write $\alpha = (\alpha^{(1)},\alpha^{(2)}) \in \Z_+^{n_1} \times \Z_+^{n_2}$.  For $m_1,m_2 \in \R$, we define the symbol class $S^{m_1,m_2}$ to be the space of all smooth functions $a(x,\xi)$ such that
\begin{equation}
 \sup_{x,\xi}|\partial_x^\beta\partial_\xi^\alpha a(x,\xi)| \leq C_{\alpha,\beta}\left<\xi^{(1)}\right>^{m_1-|\alpha^{(1)}|}\left<\xi^{(2)}\right>^{m_2-|\alpha^{(2)}|}.
\end{equation}
The space $S^{m_1,m_2}$ is a Fr\'echet space whose topology is generated by the seminorms
\begin{equation}
 \|a\|_{S^{m_1,m_2}_{\alpha,\beta}} := \sup_{x,\xi}\left<\xi^{(1)}\right>^{-m_1+|\alpha^{(1)}|}\left<\xi^{(2)}\right>^{-m_2+|\alpha^{(2)}|}|\partial_x^\beta\partial_\xi^\alpha a(x,\xi)|
\end{equation}
We define $OS^{m_1,m_2} = OS^{m_1,m_2}(\R^{n_1}\times\R^{n_2})$ to be the class of all operators obtained via the formula (\ref{quantize}) for $a \in S^{m_1,m_2}$.  An operator in $OS^{m_1,m_2}$ is called a \textit{product type pseudodifferential operator}.

For the purposes of generalizing Theorem \ref{ThmPSDO} to our anisotropic spaces, we will need to introduce yet another type of symbol class. These symbols are ``anisotropic symbols", since they obey an anisotropic type decay.  Namely, given $m_1,m_2 \in \R$, define the symbol class $S^{(m_1,m_2)}$ to be the space of all smooth functions $a(x,\xi)$ such that
\begin{equation}
 \sup_{x,\xi}|\partial_x^\beta\partial_\xi^\alpha a(x,\xi)| \leq C_{\alpha,\beta}\left<\xi\right>^{m_1-|\alpha^{(1)}|}\left<\xi^{(2)}\right>^{m_2-|\alpha^{(2)}|}.
\end{equation}
We define the seminorms $\|\cdot\|_{S^{(m_1,m_2)}_{\alpha,\beta}}$ on $S^{(m_1,m_2)}$ in the obvious way.  Thus, when we differentiate symbols in $S^{(m_1,m_2)}$ in the $\xi^{(1)}$ variables, we get full radial decay in $\xi$, but we only get decay in $\xi^{(2)}$ when we differentiate in the $\xi^{(2)}$ derivatives.  Hence, we have the containments, $S^0 \subset S^{(0,0)} \subset S^{0,0}$, where the symbol classes obey radial, anisotropic, and product type decay upon differentiation in the $\xi$ variables, respectively.   Define the class of anisotropic type operators $OS^{(m_1,m_2)} = OS^{(m_1,m_2)}(\R^{n_1}\times\R^{n_2})$ in the obvious way.  In fact, all the operators in this paper will be of anisotropic type; we only consider them as product type, when applicable, in order to make use of the mapping properties of product type operators as established in \cite{Ya2}.  For all $m \in \R$, we have the obvious inclusions
\begin{align*}
  OS^m \subset OS^{(m,0)} \subset OS^{m,m}.
\end{align*}

We now have the following theorem:

\begin{Theorem}\label{ThmProdPSDO} \hspace{.1in}
\begin{enumerate}
 \item For all $m_1,m_1',m_2,m_2' \in \R$, we have the composition rules
 \begin{align*}
   OS^{m_1,m_2}\circ OS^{m_1',m_2'} \to OS^{m_1+m_1',m_2+m_2'}.\\
   OS^{(m_1,m_2)}\circ OS^{(m_1',m_2')} \to OS^{(m_1+m_1',m_2+m_2')}.
 \end{align*}
 \item If $P \in OS^{0,0}$ then $P$ is a bounded operator on $A^{(s_1,s_2),p}$ for all $s_1,s_2 \in \R$ and $1<p<\infty$.  Moreover, for any fixed $s_1$,$s_2$,$p$, the operator norm of $P$ is bounded in terms of only finitely many symbol semi-norms $S^{0,0}_{\alpha,\beta}$.
 \item If $P \in OS^{(m_1,m_2)}$, then $P: A^{(s_1,s_2),p} \to A^{(s_1-m_1,s_2-m_2),p}$ is bounded.  Moreover, the norm of $P$ depends on only finitely many semi-norms $S^{(m_1,m_2)}_{\alpha,\beta}$.
\end{enumerate}
\end{Theorem}

\begin{proof}(i) The second composition rule is exactly \cite[Proposition 14.32]{WRL}. One checks that the proof given there also follows through verbatim for the product type operators in $OS^{m_1,m_2}$.

(ii) By definition of $A^{(s_1,s_2),p}$, the operator $P$ is bounded on $A^{(s_1,s_2),p}$ if and only if $\tilde P:= J_{(2)}^{-s_2}PJ_{(2)}^{s_2}$ is bounded on $A^{s_1,p}$.  Since $J_{(2)}^{\pm s_2} \in OS^{(0, \pm s_2)}$, then $\tilde P \in OS^{(0,0)}$ by (i).  By \cite{Ya2}, elements of $OS^{0,0}$ are bounded on product type spaces.  In particular, the operator $\tilde P$ is bounded on the product spaces $F^{k,0}_{p,2}$ and $F^{0,k}_{p,2}$ for all nonnegative integers $k$.  Thus, $\tilde P$ is bounded on
$$F^{k,0}_{p,2} \cap F^{0,k}_{p,2} = H^{k,p}, \qquad k \geq 0$$
where the above equality is a simple consequence of Lemma \ref{LemmaProdLP}.  By duality, $P$ is bounded on $H^{k,p}$ for all $k \in \Z$.  Since all the $A^{s,p}$ spaces arise as interpolation spaces from the $H^{k,p}$ for $k$ an integer (see \cite{Tr}), it follows that $\tilde P$ is bounded on $A^{s_1,p}$.  Thus, $P$ is bounded on $A^{(s_1,s_2),p}$.  For the final statement, \cite{Ya2} shows that the norm of a product-type pseudodifferential operator on product-type spaces depend only on finitely many of the semi-norms of the symbol.  Moreover, when composing symbols, if one inspects the proof of (i), one sees that each symbol semi-norm of a composite symbol depends continuously on only finitely many semi-norms of the symbols of each factor.  The statement now follows.

(iii) The proof is similar to (ii).  Namely, if we define
$$\tilde P = J^{(s_1-m_1)}J^{(s_2-m_2)}_{(2)}PJ^{-s_1}J_{(2)}^{-s_2},$$
then $\tilde P \in OS^{(0,0)}$ by the composition rule for anisotropic type pseudodifferential operators in (i).  It now suffices to show that $\tilde P$ is bounded on $A^{0,p}$.  But this follows from $OS^{(0,0)} \subset OS^{0,0}$ and (ii).
\end{proof}

\section{Basic Properties}

In addition to the pseudodifferential properties of the previous section, we wish to establish some other basic properties for our anisotropic function spaces.  These properties are modeled off the most basic properties that hold for their isotropic counterparts in \cite{Tr1} and \cite{Tr}.
Of course, since the isotropic spaces are very well understood, we will carry out our proofs with a minimal amount of effort by reducing the situation to the isotropic case as much as possible.  Using the pseudodifferential mapping properties in the previous section, this is often done very easily.

\begin{Lemma}\label{LemmaBasic} Let $s_1,s_2 \in \R$ and $1<p<\infty$.
  \begin{enumerate}
    \item Smooth functions are dense in $A^{(s_1,s_2),p}(\R^{n_1}\times\R^{n_2})$.
    \item Let $\varphi \in C^\infty_0(\R^n)$.  Then multiplication by $\varphi$ defines a bounded operator on $A^{(s_1,s_2),p}(\R^{n_1}\times\R^{n_2})$.
    \item Let $\Phi_j: \R^{n_j} \to \R^{n_j}$ be diffeomorphisms such that $D^\alpha\Phi_j \in L^\infty$ for all $|\alpha| \geq 1$ and $\inf_{x \in \R^{n_j}}|\det D\Phi_j(x)| \geq c$ for some $c > 0$, $j=1,2$.  Then $f \mapsto f \circ (\Phi_1\times\Phi_2)$ defines a bounded operator on $A^{(s_1,s_2),p}(\R^{n_1}\times \R^{n_2})$.
  \end{enumerate}
\end{Lemma}

\begin{proof}
(i) Consider the following approximations to the identity, $\varphi_\eps(x) = \eps^{-n}\varphi(x/\eps)$, $\eps > 0$, where $\varphi \in C^\infty_0(\R^n)$ integrates to unity.  For any $f \in A^{(s_1,s_2),p}$, let $f_\eps = \varphi_\eps*f$ denote the convolution with $\varphi_\eps$. Note that this convolution operation is also given by a zeroth order pseudodifferential operator with total symbol $\hat \varphi(\eps \xi)$.  Since $\hat\varphi(\eps \xi) \to 1$ as a zeroth order symbol, by Theorem \ref{ThmProdPSDO}, the result follows, since then the $f_\eps$ are smooth and converge to $f$ in $A^{(s_1,s_2),p}$ as $\eps \to 0$.

(ii) Multiplication by a smooth function is also given by a pseudodifferential operator and so Theorem \ref{ThmProdPSDO} establishes the result again.

(iii) To show that $(\Phi_1 \times \Phi_2)^*: A^{(s_1,s_2),p} \to A^{(s_1,s_2),p}$ is bounded, it suffices to show that
\begin{equation}
 J_{(2)}^{s_2}(\Phi_1 \times \Phi_2)^*J_{(2)}^{-s_2}: A^{s_1,p} \to A^{s_1,p} \label{JDJ}
\end{equation}
is bounded.  However, since a pseudodifferential operator transforms under a diffeomorphism satisfying the above hypotheses to another pseudodifferential operator of the same order (see \cite[Theorem 7.16]{WRL}), we see that the map (\ref{JDJ}) is given by the composite map $J_{(2)}^{s_2}\tilde J_{(2)} (\Phi_1 \times \Phi_2)^*$, where $\tilde J_{(2)}$ is the conjugate of $J_{(2)}^{-s_2}$ by the diffeomorphism $\Phi_2$.  Thus, since $\tilde J_{(2)} \in OS^{0,-s_2}$, we have $J_{(2)}^{s_2}\tilde J_{(2)} \in OS^{0,0}$ by Theorem \ref{ThmProdPSDO}.  The map $(\Phi_1 \times \Phi_2)^*$ is bounded on the isotropic space $A^{s_1,p}$ by \cite[Proposition 4.3.1]{Tr2}, and so now the boundedness of (\ref{JDJ}) follows.
\end{proof}

Before proceeding to further properties of our anisotropic function spaces, we will provide some simple equivalent characterizations of our function spaces in terms of other standard function spaces.  For this, we recall some definitions.  Given an arbitrary Banach space $\X$, recall that we can define vector valued Sobolev and Bessel potential spaces:
\begin{align*}
  L^p(\R^n, \X) & := \{f \in \mathcal{S}'(\R^n,\X) : \|f\|_{L^p(\R^n,\X)} = \big\| \|f(\cdot)\|_{\X}^p \big\|_{L^p(\R^n)} < \infty\} \\
  H^{s,p}(\R^n, \X) & := \{f \in \mathcal{S}'(\R^n,\X) : \|f\|_{H^{s,p}(\R^n,\X)} = \|\F^{-1}J^s\F f\|_{L^p(\R^n,\X)} < \infty\}
\end{align*}
for $1 < p < \infty$ and $s \in \R$. Here, $\mathcal{S}'(\R^n,\X)$ denotes the space of $\X$-valued tempered distributions, i.e. the space of continuous linear maps $f: \mathcal{S}(\R^n) \to \X$.  One can also define vector-valued Besov spaces and Triebel-Lizorkin spaces in a similar fashion.  The study of vector valued function spaces has an entire literature of its own, as new difficulties arise concerning operator-valued multipliers and to what extent properties of classical scalar-valued function spaces continue to hold.  For some results and further reading, see e.g., \cite{GW, SchSi}.

Next, recall that for $1 < p < \infty$, there is a $p$-nuclear tensor norm $\alpha_p$, which is a uniform cross norm one can define on the tensor product of any two Banach spaces $\X$ and $\Y$.  Namely, if we let $z \in \X \otimes \Y$ be given by
\begin{equation}
  z = \sum_{j=1}^n x_i\otimes y_i, \qquad x_i \in \X, y_i \in \Y, \label{repz}
\end{equation}
then we define
$$\alpha_p(z,\X,\Y) = \inf \left\{ \left(\sum_{i=1}^n\|x_i\|_\X^p\right)^{1/p} \cdot \sup \left\{ \left( \sum_{i=1}^n \|\psi(y_i)\|_{\Y}^{p'} \right)^{1/p'} : \psi \in \Y',\; \|\psi\|_{\Y'} \leq 1\right\} \right\}$$
where $\frac{1}{p} + \frac{1}{p'} = 1$ and the infimum is over all representations of $z$ as in (\ref{repz}).

By using the so called Fubini property of the isotropic Bessel potential and Besov spaces, and the results on tensor products of Banach spaces in \cite{SiUl}, we have the following:

\begin{Lemma}\label{LemmaEquiv}
  Let $1 < p < \infty$.
  \begin{enumerate}
    \item If $s_1,s_2 \geq 0$ then
    \begin{align}
      H^{(s_1,s_2),p}(\R^{n_1}\times\R^{n_2}) & = L^p(\R^{n_1}, H^{s_1+s_2,p}(\R^{n_2})) \cap H^{s_2,p}(\R^{n_1}, H^{s_1,p}(\R^{n_2})) \label{equiv1} \\
      & = \bigcap_{0 \leq s \leq s_1} H^{s_1 - s,p}(\R^{n_1}, H^{s_1+s_2,p}(\R^{n_2})) \label{equiv2} \\
      &= \bigcap_{0 \leq s \leq s_1} H^{s_1 - s,p}(\R^{n_1}) \otimes_{\alpha_p} H^{s_2+s,p}(\R^{n_2}). \label{equiv3}
    \end{align}
For $p = 2$, the tensor product in (\ref{equiv3}) can be replaced with just the ordinary tensor product of Hilbert spaces.
\item If $s_1 > 0$ and $s_2 \geq 0$, we have
\begin{align}
  B^{(s_1,s_2),p}(\R^{n_1}\times\R^{n_2}) & = L^p(\R^{n_1}, B^{s_1+s_2,p}(\R^{n_2})) \cap H^{s_2,p}(\R^{n_2}, B^{s_1,p}(\R^{n_1})).
\end{align}
  \end{enumerate}
\end{Lemma}

\begin{proof} (i) For $s_1 \geq 0$, we have
$$H^{s_1,p}(\R^{n_1}\times\R^{n_2}) = L^p(\R^{n_1}, H^{s_1,p}(\R^{n_2})) \cap L^p(\R^{n_2}, H^{s_1,p}(\R^{n_1})).$$
This follows from the lift property Lemma \ref{LemmaLift} and fact that the multiplier $\frac{\left<\xi\right>^{s_1}}{\left<\xi_{(1)}\right>^{s_1} + \left<\xi_{(2)}\right>^{s_1}}$ and its reciprocal are both elements of $OS^{0,0}$.  We now apply $J_{(2)}^{-s_2}$ to both sides of the above to obtain (\ref{equiv1}).  The next equality follows from intersecting (\ref{equiv1}) with $H^{s_1-s,p}(\R^{n_1}, H^{s_2+s,p}(\R^{n_2})$, for $0 < s < s_1$.  We have
\begin{equation}
  H^{(s_1,s_2),p}(\R^{n_1}\times\R^{n_2}) \subseteq H^{s_1-s,p}(\R^{n_1}, H^{s_2+s,p}(\R^{n_2})) \label{eq:contain}
\end{equation}
since $\frac{\left<\xi_{(1)}\right>^{s_1-s}\left<\xi_{(2)}\right>^{s_2+s}}{\left<\xi\right>^{s_1}\left<\xi_{(2)}\right>^{s_2}} \in OS^{0,0}$, and so (\ref{equiv1}) implies (\ref{equiv2}).  Finally, (\ref{equiv3}) follows from the fact that
\begin{equation}
  H^{s,p}(\R^{n_1}, H^{t,p}(\R^{n_2})) = F^{s,t}_{p,2}(\R^{n_1} \times \R^{n_2}) = H^{s,p}(\R^{n_1}) \otimes_{\alpha_p} H^{t,p}(\R^{n_2}), \qquad s,t \in \R. \label{Hst}
\end{equation}
Here, the first equality follows straight from the definitions and the second equality follows from \cite[Proposition 3.1]{SiUl}.  For $p = 2$, the final statement follows easily from the unitarity of the Fourier transform and the fact that $L^2(\R^{n_1} \times \R^{n_2}) = L^2(\R^{n_1}) \otimes L^2(\R^{n_2})$, where $\otimes$ is the ordinary tensor product of Hilbert spaces.

(ii) For $s_1 > 0$, the Besov space $B^{s_1,p}(\R^{n_1} \times \R^{n_2})$ satisfies the following Fubini property (see \cite[Theorem 2.5.13]{Tr}):
$$B^{s_1,p}(\R^{n_1} \times \R^{n_2}) = L^p(\R^{n_1}, B^{s_1,p}(\R^{n_2})) \cap L^p(\R^{n_2}, B^{s_1,p}(\R^{n_1})).$$
We now apply $J_{(2)}^{-s_2}$ to both sides.
\end{proof}

Having defined our anisotropic function spaces on products of Euclidean space, we can now define them on products of subsets of Euclidean space in the usual way.  Namely, given subsets $\Omega_i \subset \R^{n_i}$, $i=1,2$, we define $A^{(s_1,s_2),p}(\Omega_1\times\Omega_2)$ to be the space of restrictions of elements of $A^{(s_1,s_2),p}(\R^{n_1}\times\R^{n_2})$ to $\Omega_1\times\Omega_2$.  In other words, given $f \in A^{(s_1,s_2),p}(\R^{n_1}\times\R^{n_2})$, it restricts to an element of $(C^\infty_0(\Omega_1\times\Omega_2))'$, the space of distributions on $\Omega_1\times\Omega_2$, in the natural way.  Call this element $r_{\Omega_1\times\Omega_2}(f)$.  Then we define
$$A^{(s_1,s_2),p}(\Omega_1\times\Omega_2) = \{f \in (C^\infty_0(\Omega_1\times\Omega_2))' : f = r_{\Omega_1\times\Omega_2}(g),\; g \in A^{(s_1,s_2),p}(\R^{n_1}\times\R^{n_2})\}$$
and we equip this space with the norm
$$\|f\|_{A^{(s_1,s_2),p}(\Omega_1\times\Omega_2)} = \inf_{\{g\;:\; r_{\Omega_1\times\Omega_2}(g)=f\}}\|g\|_{A^{(s_1,s_2),p}(\R^{n_1}\times\R^{n_2})}.$$

The next theorem works out the trace and extension theorems for anisotropic spaces, where one wishes to understand how restricting a function on $\R^n$ to a hyperplane $\R^{n-1}$ behaves and likewise for extending a funtion on $\R^{n-1}$ to $\R^n$.  There are two cases to consider.  The first and trivial case is when the anisotropy is tangential to the boundary.  Then the anisotropy simply commutes with the trace map to the boundary.  The nontrivial case is when there is some anisotropy in the direction normal to the hyperplane. In this case, the below theorem tells us if there is enough anisotropy, then taking a trace ``costs'' derivatives in mostly the anisotropic directions; in fact if $p \geq 2$, only the regularity in the anisotropic directions are decreased.

In the following, let $\R^n_+$ denote the half space $\{(x_1,\ldots, x_n) \in \R^n: x_n > 0\}$.  Given a function $f$ defined on $\R^n_+$ and smooth up to the boundary, we can define its $m$th order trace map on $\R^{n-1} = \partial \R^n_+$ via
\begin{equation}
 r_m: f \mapsto (f, \partial_{x_n}f, \ldots, \partial_{x_n}^mf)|_{\R^{n-1}}.
\end{equation}
Likewise, we have trace maps from $\R^{n_1}_+ \times \R^{n_2}$ and $\R^{n_1} \times \R^{n_2}_+$ to $\partial \R^{n_1}_+\times\R^{n_2}$ and $\R^{n_1}\times\partial\R^{n_2}_+$, respectively.  These trace maps on spaces of smooth functions extend to anisotropic function spaces in the following way:

\begin{Theorem}\label{ThmATrace}(Anisotropic Traces and Extensions)
   \begin{enumerate}
     \item (Tangential anisotropy) For $s_1 > m + 1/p$ and $s_2 \geq 0$, the $m$th order trace map from $\R^{n_1}_+\times\R^{n_2}$ to $\partial\R^{n_1}\times\R^{n_2}$ satisfies
     $$r_m: A^{(s_1,s_2),p}(\R^{n_1}_+\times\R^{n_2}) \to \oplus_{j=0}^{m} B^{(s_1-j-1/p,s_2),p}(\partial\R^{n_1}\times\R^{n_2}).$$
     Furthermore, for all $s_1 \in \R$, there exists a boundary extension map
     $$e_m: \oplus_{j=0}^{m} B^{(s_1-j-1/p,s_2),p}(\partial\R^{n_1}\times \R^{n_2}) \to A^{(s_1,s_2),p}(\R^{n_1}_+\times \R^{n_2}),$$
     and for $s_1 > m + 1/p$, we have $r_me_m = \mathrm{id}$.  Moreover, for every $k \in \N$, we have an extension map
   \begin{equation}
     E_k: A^{(s_1,s_2),p}(\R^{n_1}_+\times \R^{n_2}) \to A^{(s_1,s_2),p}(\R^{n_1}\times \R^{n_2}), \label{Ek}
   \end{equation}
   for $|s_1| < k$.
     \item (Mixed anisotropy) For $s_1 \geq 0$, $s_2 >m+1/p$, the $m$th order trace map from $\R^{n_1}\times\R^{n_2}_+$ to $\R^{n_1}\times\partial\R^{n_2}_+$ satisfies
     \begin{align*}
     	r_m: H^{(s_1,s_2),p}(\R^{n_1}\times\R^{n_2}_+) & \to \oplus_{j=0}^{m-1}H^{(s_1, s_2-j-1/p-\eps_2),p}(\R^{n_1}\times\partial\R^{n_2}_+)\\
	r_m: B^{(s_1,s_2),p}(\R^{n_1}\times\R^{n_2}_+) & \to \oplus_{j=0}^{m-1}B^{(s_1-\eps_1, s_2-j-1/p-\eps_2),p}(\R^{n_1}\times\partial\R^{n_2}_+),
     \end{align*}
where $\eps_1$ and $\eps_2$ satisfy the following:
\begin{enumerate}
\item if $p> 2$, then $\eps_1=0$ and $\eps_2 > 0$ is arbitrary;
\item if $p=2$, then $\eps_1=\eps_2=0$;
\item if $1<p<2$, then $\eps_1,\eps_2>0$ are arbitrary.
   \end{enumerate}
\end{enumerate}
\end{Theorem}

\begin{proof}
(i) The isotropic case $s_2 = 0$ is well documented, e.g., see \cite[Theorem 2.7.2]{Tr1}.  For the anisotropic case $s_2 > 0$, the result for $r_m$ follows from the isotropic case because (pseudodifferential) operators on $\R^{n_2}$,  which thus act tangentially to the boundary, commute with $r_m$.  For the extension maps $e_m$ and $E_m$, one can also see from their constructions in \cite{Tr1,Tr} that tangential pseudodifferential operators commute with them, and so the anisotropic case also follows from the isotropic case.

(ii) Unlike the previous case, here some work needs to be done. For this, we make use of the product type spaces considered in Definition \ref{DefProdSpaces}.  This is because then we can make use of the trace theorems found in \cite{ST} for these spaces.  To facilitate this, for all integers $n_1,n_2 \geq 1$, let us define
\begin{align*}
SH^{s_1,s_2,p}(\R^{n_1}\times\R^{n_2}) &:= F^{s_1,s_2}_{p,2}(\R^{n_1}\times\R^{n_2})\\
SB^{s_1,s_2,p}(\R^{n_1}\times\R^{n_2}) &:= B^{s_1,s_2}_{p,p}(\R^{n_1}\times\R^{n_2}),
\end{align*}
where the notation is a hybrid between the notation for our anisotropic spaces and the notation for spaces of dominating mixed smoothness in \cite{ST}. Without loss of generality, let $m = 0$ and let $r = r_0$ be the zeroth order trace map. Applying \cite[Theorem 2.4.2]{ST}\footnote{In \cite{ST}, the results are stated for $n_1=n_2=1$, but they easily generalize to all $n_1,n_2 \geq 1$.}, we have the following bounded trace maps for spaces of dominating mixed smoothness
\begin{align}
 r: SH^{s_1,s_2,p}((\R^{n_1}\times \partial\R^{n_2}_+) \times \R_+) \to H^{s_1,p}(\R^{n_1}\times\partial\R^{n_2}_+) \label{rSH}\\
 r: SB^{s_1,s_2,p}((\R^{n_1}\times \partial\R^{n_2}_+) \times \R_+) \to B^{s_1,p}(\R^{n_1}\times\partial\R^{n_2}_+), \label{rSB}
\end{align}
whenever $s_2 > 1/p$.  Here, we regard $\R^{n_2}_+ = \partial \R^{n_2}_+\times\R_+$. Observe that the right-hand side of the above are the usual isotropic spaces on $\R^{n_1}\times\partial\R^{n_2}_+$.  Moreover, the trace maps (\ref{rSH}) and (\ref{rSB}) merely require enough regularity in the normal direction $\R_+$ to the boundary and these maps preserve the tangential regularity.

For the anisotropic spaces, we will derive our trace theorem from embedding them into the appropriate space of dominating mixed smoothness.  For the Bessel potential spaces, this is straightforward in light of tensor product representation in Lemma \ref{LemmaEquiv} and equations (\ref{eq:contain}) and (\ref{Hst}).  Namely, we write
\begin{align}
  H^{(s_1,s_2),p}(\R^{n_1}\times\R^{n_2}_+) &= \bigcap_{0 \leq s \leq s_1} H^{s_1-s,p}(\R^{n_1}) \otimes_{\alpha_p} H^{s_2+s,p}(\R^{n_2}_+) \label{Hs1s2-1}\\
  & \subseteq \bigcap_{0 \leq s \leq s_1} H^{s_1-s,p}(\R^{n_1}) \otimes_{\alpha_p} SH^{s_2+s-1/p-\eps,1/p+\eps,p}(\partial \R^{n_2}_+ \times \R_+), \label{Hs1s2-2}
\end{align}
where $\eps > 0$ is arbitrary.  Thus, when we apply the trace map $r$ to (\ref{Hs1s2-2}), equation (\ref{rSH}) and Lemma \ref{LemmaEquiv} implies that
\begin{align*}
  r: H^{(s_1,s_2),p}(\R^{n_1}\times\R^{n_2}_+) & \to  \bigcap_{0 \leq s \leq s_1} H^{s_1-s,p}(\R^{n_1}) \otimes_{\alpha_p} H^{s_2-1/p-\eps+s,p}(\partial \R^{n_2}_+ \times \R_+)\\
  & = H^{(s_1,s_2-1/p-\eps),p}(\R^{n_1} \times \partial \R^{n_2}_+).
\end{align*}
For $p = 2$, we apply the trace map $r$ to (\ref{Hs1s2-1}) directly and apply the standard isotropic trace theorem to obtain
\begin{align*}
  r: H^{(s_1,s_2),2}(\R^{n_1}\times\R^{n_2}_+) & \to  \bigcap_{0 \leq s \leq s_1} H^{s_1-s,2}(\R^{n_1}) \otimes H^{s_2+s-1/2,2}(\partial \R^{n_2}_+)\\
  & = H^{(s_1,s_2-1/2),2}(\R^{n_1}\times\partial\R^{n_2}_+).
\end{align*}

For the Besov case, we proceed as follows.  First, we have the following embeddings:
\begin{equation}
\begin{array}{rcll}
 B^{0,p} & \subseteq & SB^{0,0,p} &\;\; 2 \leq p<\infty\\
 B^{0,p} & \subseteq & SB^{-\eps,-\eps,p} &\;\; 1<p<2, \eps > 0.
\end{array} \label{embeddings}
\end{equation}
These two embeddings follow from the following computation.  Given a function $f$ on $\R^n = \R^{n_1}\times\R^{n_2}$, we have $f = \sum_k f_k$, where the $f_k = \mathcal{F}^{-1}\varphi_k\mathcal{F}f$ yield a radial Littlewood-Paley decomposition on $\R^n$, and we also have $f = \sum_{i,j}f_{i,j}$, where the $f_{i,j} = \F\varphi_i^{(1)}\varphi_j^{(2)}\F^{-1}f$ yield a product Littlewood-Paley decomposition for $f$.  We have
\begin{align}
  \|f\|_{B^{0,p}}^p &= \sum_k\|\{f_k\}_{l^p_k}\|_{L^p}^p \nonumber \\
&= \sum_k \|f_k\|_{L^p}^p \nonumber \\
  &= \sum_k \|\sum_{i,j}\{(f_k)_{i,j}\}_{\ell^2_{i,j}}\|_{L^p}^p. \nonumber
\end{align}
The last line is just the product Littlewood-Paley decomposition, given by Lemma \ref{LemmaProdLP}, applied to each $f_k$. When $p \geq 2$, then we have
\begin{align*}
\sum_k \|\sum_{i,j}\{(f_k)_{i,j}\}_{\ell^2_{i,j}}\|_{L^p}^p  & \geq  \sum_k \|\sum_{i,j}\{(f_k)_{i,j}\}_{\ell^p_{i,j}}\|_{L^p}^p\\
  & = \sum_{i,j}\sum_{k \atop |\max(i,j) -k| \leq 2}\|(f_{i,j})_k\|_{L^p}^p\\
  & \sim \sum_{i,j}\|f_{i,j}\|_{L^p}^p\\
  & = \|f\|_{SB^{(0,0),p}}^p.
\end{align*}
The first line follows from the inclusion $\ell^p \hookrightarrow \ell^2$ for $p\geq 2$, and the line after follows from a change of summation and the fact that $f_{i,j}$ has frequency support only on $\xi \sim 2^{\max(i,j)}$.  This proves the second embedding of (\ref{embeddings}).  To obtain the last embedding of (\ref{embeddings}), we instead use the inequality $\{(f_k)_{i,j}\}_{\ell^2_{i,j}} \geq C\{2^{-\eps i}2^{-\eps j}(f_k)_{i,j}\}_{\ell^p_{i,j}}$ in the first line in the above, which follows from H\"older's inequality.

So let $f \in B^{(s_1,s_2),p}(\R^{n_1}\times\R^{n_2}_+)$.  Consider first the case $p > 2$.  Let $J_{\partial}$ be the Bessel potential operator of order one on $\partial \R^{n_2}_+$.  It follows that $J_{\partial}^{s_2'}f \in B^{(s_1,s_2-s_2'),p}(\R^{n_1}\times\R^{n_2}_+)$ for $s_2' \geq 0$. Furthermore, we have
\begin{equation}
B^{(s_1,s_2-s_2'),p}(\R^{n_1}\times\R^{n_2}_+) \subseteq SB^{s_1,s_2-s_2',p}((\R^{n_1}\times\partial\R^{n_2}_+)\times\R_+), \label{HintoSH}
\end{equation}
which one can see as follows. On $\R^{n_1}\times\R^{n_2}$, we have
\begin{equation}
  B^{(0,0),p}(\R^{n_1}\times\R^{n_2}) \subseteq SB^{0,0,p}((\R^{n_1}\times \R^{n_2-1})\times\R) \label{H00=SH00}
\end{equation}
by (\ref{embeddings}).  With respect to the decomposition $\R^n = \R^{n_1}\times\R^{n_2-1}\times\R$, let $J_{n_1+n_2-1}$ and $J_1$ denote the Bessel potential operators of order one on $\R^{n_1+n_2-1}$ and $\R$, respectively. Then for $s_2' \leq s_2$, consider the isomorphisms
\begin{align*}
J^{-s_1}J_{(2)}^{-(s_2-s_2')}: B^{(0,0),p}(\R^{n_1}\times\R^{n_2}) & \cong B^{(s_1,s_2-s_2'),p}(\R^{n_1}\times\R^{n_2}) \\
 J_{n_1+n_2-1}^{-s_1}J_{1}^{-(s_2-s_2')}: SB^{0,0,p}((\R^{n_1}\times\R^{n_2-1})\times\R) & \cong SB^{s_1,s_2-s_2',p}((\R^{n_1}\times\R^{n_2-1})\times\R).
\end{align*}
These isomorphisms and (\ref{H00=SH00}) imply that the embedding (\ref{HintoSH}) follows from the operator $$J_{n_1+n_2-1}^{s_1}J_{1}^{s_2-s_2'}J^{-s_1}J_{(2)}^{-(s_2-s_2')}: B^{(0,0),p}(\R^{n_1}\times\R^{n_2}) \to B^{(0,0),p}(\R^{n_1}\times\R^{n_2})$$
being bounded.  However, the above operator is equal to
\begin{equation}
   J^{-s_1}J_{n_1+n_2-1}^{s_1} \cdot J_{(2)}^{-(s_2-s_2')}J_{1}^{s_2-s_2'}, \label{JJ}
\end{equation}
which is the product of two operators in $OS^{0,0}((\R^{n_1}\times\R^{n_2-1})\times\R)$ and $OS^{0,0}(\R^{n_2-1}\times\R)$, respectively, since $s_1, s_2 - s_2' \geq 0$. It follows from Theorem \ref{ThmProdPSDO} that (\ref{JJ}) is bounded on $B^{(0,0),p}(\R^{n_1}\times\R^{n_2})$ and hence (\ref{HintoSH}) holds.

Thus, we have $J_{\partial}^{s_2'}f \in SB^{s_1,s_2-s_2',p}((\R^{n_1}\times\partial\R^{n_2}_+)\times\R_+)$.  So if $s_2 - s_2' > 1/p$, then $$J_{\partial}^{s_2'}r(f) = r(J_{\partial}^{s_2'}f) \in B^{s_1,p}(\R^{n_1}\times\partial\R^{n_2}_+),$$
by (\ref{rSH}). This means
$$r(f) \in B^{(s_1,s_2'),p}(\R^{n_1}\times\partial \R^{n_2}_+)$$
for $s_2' < s_2-1/p$. This proves the trace theorem for the anisotropic Bessel potential spaces when $2<p<\infty$.  Similarly, for $1<p<2$, the same argument above shows that (\ref{embeddings}) implies
$$B^{(s_1,s_2-s_2'),p}(\R^{n_1}\times\R^{n_2}_+) \subseteq SB^{s_1-\eps,s_2-s_2'-\eps,p}((\R^{n_1}\times\partial\R^{n_2}_+)\times\R_+)  \qquad \eps > 0.$$
In the same way, we obtain the conclusion of the theorem for anisotropic Besov spaces.
\end{proof}

Next, we want to understand the various interpolation properties of these function spaces in the regularity parameters.  Recall that for $0 < \theta < 1$, there is a complex interpolation functor $[\cdot,\cdot]_\theta$ mapping the category of interpolation of couples of Banach spaces to Banach spaces.  For isotropic spaces, we have $[A^{s_0,p},A^{s_1,p}]_\theta = A^{(1-\theta)s_0+\theta s_1,p}$ for all $s_0,s_1 \in \R$.  We have the following anisotropic generalization:

\begin{Lemma}\label{LemmaInterp}
  Let $s_1,s_1' \in \R$ and $s_2,s_2' \in \R$, and $0 < \theta < 1$.  Furthermore, suppose $s_i \leq s_i'$, $i=1,2$.  Then we have
  $$[A^{(s_1,s_2),p},A^{(s_1',s_2'),p}]_\theta = A^{((1-\theta)s_1+\theta s_1',(1-\theta)s_2+\theta s_2'),p}.$$
\end{Lemma}

\begin{proof}
We can think of $A^{(s_1',s_2'),p} \subset A^{(s_1,s_2),p}$ as the domain of the unbounded operator $\Lambda = J^{s_1'-s_1}J_{(2)}^{s_2'-s_2}$.  Using Theorem \ref{ThmProdPSDO}, one can verify that $\Lambda$ satisfies all the hypotheses of \cite[Theorem 1.15.3]{Tr}.  We now conclude from that theorem that $[A^{(s_1,s_2),p},A^{(s_1',s_2'),p}]_\theta$ is precisely the domain of the operator $\Lambda^\theta$, which is precisely space $A^{((1-\theta)s_1+\theta s_1',(1-\theta)s_2+\theta s_2'),p}$ by the lift property.
\end{proof}

\section{Multiplication}

There are numerous multiplication theorems concerning the usual isotropic Bessel potential and Besov spaces, see e.g. \cite{RS}.  In passing to their anisotropic counterparts, we want to see how much of the anisotropy can be preserved after multiplication.  When $A^{(s_1,s_2),p}(\R^n)$ is such that $A^{s_1,p}(\R^n)$ is a Banach algebra, then with very little work we have that the anisotropy $s_2$ is completely preserved.

\begin{Lemma}
 Let $s_1 > n/p$.  Then for all $s_2 \geq 0$, we have that $A^{(s_1,s_2),p}(\R^{n_1}\times\R^{n_2})$ is a Banach algebra.  In particular, for $s_2'$,$s_2'' \geq 0$ we have a map
\begin{align*}
A^{(s_1,s_2'),p}(\R^{n_1}\times\R^{n_2}) \times A^{(s_1,s_2''),p}(\R^{n_1}\times\R^{n_2}) \to A^{(s_1,\min(s_2',s_2'')),p}(\R^{n_1}\times\R^{n_2}).
\end{align*}
\end{Lemma}

\begin{proof}
 For $s_2 = 0$, this result is standard.  For $s_2$ a nonnegative integer, this now follows from the ordinary Leibnitz rule for differentiation.  For all other $s_2 > 0$, we can interpolate between the bilinear forms (see \cite[Chapter 1.19.5]{Tr})
\begin{align}
A^{(s_1,0),p}(\R^n) \times A^{(s_1,0),p}(\R^n) & \to A^{(s_1,0),p}(\R^n) \\
A^{(s_1,k),p}(\R^n) \times A^{(s_1,k),p}(\R^n) & \to A^{(s_1,k),p}(\R^n)
\end{align}
given by multiplication, where $k > s_2$ is an integer. Using Lemma \ref{LemmaInterp}, this shows that $A^{(s_1,s_2),p}(\R^n)$ is also an algebra.
\end{proof}

When $0 < s_1 < n/p$, one has to work harder.  In the isotropic case, there are two typical multiplication maps one considers.  Namely, one can map into a space of lower regularity,
\begin{align}
 A^{s_1,p}(\R^n) \times A^{s_1,p}(\R^n) \to A^{2s_1-n/p,p}(\R^n),
\end{align}
or one can consider only those distributions that lie in $L^\infty(\R^n)$, in which case,  the algebra property is restored (see e.g. \cite[Chapter 2]{Tay}):
\begin{align}
 (A^{s_1,p}(\R^n) \cap L^\infty) \times (A^{s_1,p}(\R^n) \cap L^\infty) \to
 (A^{s_1,p}(\R^n) \cap L^\infty). \label{algebra}
\end{align}
The proofs involved in all the above multiplication theorems involve the paraproduct calculus, whereby one separates a product of functions into different frequency pieces and estimates each of these pieces appropriately.  If one tries to repeat this procedure in the anisotropic case, this procedure needs to be redone carefully.  This is because in the isotropic case, the multiplier $\left<\xi\right>$ used to define the isotropic Bessel potential and Besov spaces is radial and the frequency decomposition involved in the paraproduct calculus is also a radial decomposition.  Thus, the dyadic radial decomposition of a function is  commensurate with the norm of a function in an isotropic space.  With anisotropy however, the weight $\left<\xi_{(2)}\right>$ used to define the norm of an anisotropic space no longer has a uniform size along a radial frequency annulus.  One therefore ought to use a product frequency decomposition when doing the paraproduct calculus for the product of two functions in anisotropic spaces.  However, this would a priori result in a multiplication theorem for product type Besov and Bessel potential spaces.  For the former, these spaces do not give rise to the usual Besov spaces (i.e. $B^{0,0}_{p,p}(\R^{n_1}\times\R^{n_2}) \neq B^0_{p,p}(\R^n)$), whereas for the latter, we do in fact have $H^{0,0}_{p,2}(\R^n) = H^0_{p,2}(\R^n) = L^p(\R^n)$ by Lemma \ref{LemmaProdLP}.

To make matters simple then, we consider the case $p=2$, in which case, the anisotropic Besov and Bessel potential spaces $H^{(s_1,s_2),2}$ and $B^{(s_1,s_2),2}$ coincide and we can use the product Littlewood-Paley decompositions with impunity.  We state one particular anisotropic multiplication below, for the sake of specificity and for its application in \cite{N2}.  Even with the restrictive range placed on the function space parameters in the proposition below, the estimates are already involved due to the anisotropy.  Nonetheless, one could use the methods here to obtain other multiplication theorems for other parameters as desired.

\begin{Proposition}
 Let $s_1 > n_1/2$ and $s_2',s_2''\geq 0$.  Let $s_2 \leq \min(s_2',s_2'')$ satisfy $s_2 < s_1+s_2'+s_2''- (n_1+n_2)/2$.  Then we have a multiplication map
\begin{align}
\begin{split}
 \left(H^{(s_1,s_2'),2}(\R^{n_1}\times\R^{n_2}) \cap L^\infty\right) \times & \left(H^{(s_1,s_2''),2}(\R^{n_1}\times\R^{n_2}) \cap L^\infty\right) \to\qquad\qquad\\ & \qquad\qquad H^{(s_1,\max(s_2,0)),2}(\R^{n_1}\times\R^{n_2}).
\end{split} \label{eq:mult}
\end{align}
\end{Proposition}

\begin{proof}
We assume $s_2 > 0$, else (\ref{algebra}) implies the result. Let us introduce the following notation for convenience.  Let $\{\varphi_i(\xi^{(1)})\}_{i=0}^\infty$ be the Littlewood-Paley partition of unity on $\R^{n_1}$ defined as in Section 2.  We have $\supp \varphi_0 \subset [0,2]$ and $\supp \varphi_i = [2^{i-1},2^{i+1}]$ for all $i\geq 1$.  Likewise, let $\{\psi_j(\xi^{(2)})\}_{j=0}^\infty$ be the corresponding Littlewood-Paley decomposition on $\R^{n_2}$, defined similarly.  Let
\begin{align*}
  \varphi_i f &= \mathcal{F}^{-1}\varphi_i(\xi^{(1)})\F f\\
  \psi_j f &= \mathcal{F}^{-1}\psi_j(\xi^{(2)})\F f
\end{align*}
denote the Littlewood-Paley components of $f$ in the $\R^{n_1}$ and $\R^{n_2}$ variables, respectively. One of the key properties of this decomposition is that for any two functions $f$ and $g$, we have
$$\supp(\varphi_i f \varphi_{i'}g) \subset [2^{i-2},2^{i-3}]\times\R^{n_2}$$
if $i' \leq i - 3$.  Indeed, this follows from the support properties of the $\varphi_i$ and the fact that the Fourier transform of multiplication is convolution.  Likewise for the $\psi_j$.  Thus, if we write
\begin{equation}
  fg = \sum_{i,i',j,j'} (\varphi_i\psi_j f)(\varphi_{i'}\psi_{j'}g), \label{freqsum}
\end{equation}
we want to estimate $fg$ by grouping those terms of the above sum according to their frequencies in the appropriate manner.  Define the operators
\begin{align*}
  \Phi_i f & = \sum_{k=0}^i\varphi_kf = \F^{-1}\varphi_0(2^{-i}\xi^{(1)})\F f\\
  \Psi_i f  & = \sum_{k=0}^i\psi_kf = \F^{-1}\psi_0(2^{-i}\xi^{(2)})\F f,
\end{align*}
which are the sum of the dyadic components of $f$ whose frequencies in the $\R^{n_1}$ and $\R^{n_2}$ directions have norm less than $\sim 2^i$, respectively.

We have the following groupings for the terms of (\ref{freqsum}).  Given two frequency indices $i$ and $i'$, define $i \gtrsim i'$, $i \sim i'$, and $i \lesssim i'$ to denote $i \geq i' - 3$, $|i - i'| < 3$, and $i \leq i' - 3$, respectively.  We may therefore consider the following cases which exhaust the possible relationships among all the frequencies:
\renewcommand{\theenumi}{\arabic{enumi}.}
\begin{enumerate}
  \item $j \gtrsim j'$
  \begin{enumerate}
    \item $i,i' \lesssim j$
    \item $i \gtrsim i'$, $i > j-3$
    \item $i \sim i'$, $i > j -3$
    \item $i \lesssim i'$, $i' > j-3$
  \end{enumerate}
  \item $j \sim j'$
  \begin{enumerate}
    \item $i \gtrsim i'$
    \item $i \sim i'$
    \item $i \lesssim i'$
  \end{enumerate}
  \item $j \lesssim j'$.
\end{enumerate}
Case 3 is symmetric with Case 1, so we will focus only on the first two cases.  All these cases are straightforward to handle except Case 1(d).  This is because in all the other cases, either $f$ or $g$ has the largest frequencies and so one may throw all the derivatives of the norm $H^{(s_1,s_2),2}$ onto the dominant term.  In Case 1(d) however, $f$ dominates in the $\R^{n_2}$ direction (since $j \gtrsim j'$), but $g$ dominates in the $\R^{n_1}$ direction (since $i' \gtrsim i$).  Let us deal with the easy cases first.

Starting with Case 1(a), we do the following.  We have
\begin{align}
  \sum_{i,i',j' \lesssim j}(\varphi_i\psi_jf)(\varphi_{i'}\psi_{j'}g) &= \sum_{j}(\Phi_{j-3}\psi_j f)(\Phi_{j-3}\Psi_{j-3}g), \label{sum1a}
\end{align}
where each term in the right-most sum has frequency support on $[0,2^{j+2}] \times [2^{j-2},2^{j+2}]$.  Thus, when we take the $H^{(s_1,s_2),2}(\R^{n_1}\times\R^{n_2})$ norm, the weight $\left<\xi\right>^{s_1}\left<\xi^{(2)}\right>^{s_2}$ is dominated by an absolute constant times $2^{js_1}2^{js_2}$ on the frequency support of the $j$th term of (\ref{sum1a}).  Hence, we have
\begin{align}
  \|\sum_{j}(\Phi_{j-3}\psi_j f)(\Phi_{j-3}\Psi_{j-3}g)\|_{H^{(s_1,s_2),2}}^2 &\leq C\sum_j 2^{2j(s_1+s_2)}\|(\Phi_{j-3}\psi_j f)(\Phi_{j-3}\Psi_{j-3}g)\|_{L^2}^2 \nonumber \\
  &\leq C\sum_j 2^{2j(s_1+s_2)}\|(\Phi_{j-3}\psi_j f)\|_{L^2}^2\|(\Phi_{j-3}\Psi_{j-3}g)\|_{L^\infty}^2 \nonumber\\
    &\leq C\|g\|_{L^\infty}^2\sum_j 2^{2j(s_1+s_2)}\|(\Phi_{j-3}\psi_j f)\|_{L^2}^2\nonumber \\
    & \leq C\|g\|_{L^\infty}^2\|f\|^2_{H^{(s_1,s_2),2}}. \label{est1a}
\end{align}
In the third line above, we used the fact that
\begin{align*}
  \|\Phi_{j-3}\Psi_{j-3}g\|_{L^\infty} &\leq \|\mathcal{F}(\Phi_{j-3}\Psi_{j-3})\|_{L^1}\|g\|_{L^\infty}\\
  & \leq C\|g\|_{L^\infty},
\end{align*}
since $\|\mathcal{F}(\Phi_{j-3}\Psi_{j-3})\|_{L^1}$ is independent of $j$.

Observe that (\ref{est1a}) yields for us the multiplication theorem for the part of $fg$ that lies in the frequency portion covered by 1(a).  It remains to handle the other cases.  As mentioned, except for Case 1(d), all of them are handled similarly, since there is always at least one function whose frequency component ``dominates".  For instance, with Case 2(b), we need to estimate the sum $\sum_{i \sim i', j \sim j'}(\varphi_i\psi_j f)(\varphi_{i'}\psi_{j'}g)$.  Each term in the sum has frequency supported on $[0,2^{i+2}] \times [0,2^{j+2}]$.  Thus, by similar reasoning as before, we have
\begin{align*}
  \|\sum_{i \sim i', j \sim j'}(\varphi_i\psi_j f)(\varphi_{i'}\psi_{j'}g)\|_{H^{(s_1,s_2),2}} \leq C\|f\|_{H^{(s_1,s_2),2}}\|g\|_{L^\infty}.
\end{align*}
For Case 1(d), we proceed as follows.  We need to estimate the sum
\begin{align*}
  \sum_{i \lesssim i', j' \lesssim j \atop i' > j -3}(\varphi_i\psi_jf)(\varphi_{i'}\psi_{j'}g) = \sum_{i' > j - 3}(\Phi_{i'-3}\psi_j f)(\varphi_{i'}\Psi_{j-3}g).
\end{align*}
The $(i',j)$ term of the above sum has frequency support in $[2^{i'-2},2^{i'+2}]\times[2^{j-2},2^{j+2}]$.  The weight $\left<\xi\right>^{s_1}\left<\xi^{(2)}\right>^{s_2}$ is thus bounded by an absolute constant times $2^{s_1i'}2^{s_2j}$ on the $(i',j)$ piece (here we used $i' > j-3$).  Moreover, when we do H\"older's inequality, the weight $2^{s_1i'}$ must pair with $g$ and the weight $2^{s_2j}$ must pair with $f$ since these are the functions with the corresponding dominant frequency terms.  Namely, we estimate as follows:
\begin{align}
  \|\sum_{i' > j - 3}(\Phi_{i'-3}\psi_j f)(\varphi_{i'}\Psi_{j-3}g)\|_{H^{(s_1,s_2),2}}^2 & \leq
   C\sum_{i' > j - 3}2^{2s_1i'}2^{2s_2j}\|(\Phi_{i'-3}\psi_j f)(\varphi_{i'}\Psi_{j-3}g)\|_{L^2}^2\\
  & \leq C\sum_{i' > j - 3} \Big\{2^{2s_2j}\|\Phi_{i'-3}\psi_jf\|_{L^\infty(\R^{n_1},L^p(\R^{n_2}))}^2 \times \nonumber\\
  & \hspace{1in} \left(2^{2s_1i'}\|\varphi_{i'}\Psi_{j-3}g\|^2_{L^2(\R^{n_1},L^{q}(\R^{n_2}))}\right)\Big\}, \label{sum1d}
\end{align}
where $p$ and $q$ satisfy $1/p+1/q=2$ and are to be determined later.  If we sum over $j$ first, then since $\|\varphi_{i'}\Psi_{j-3}g\|^2_{L^2(\R^{n_1},L^{q}(\R^{n_2}))} \leq C\|\varphi_{i'}g\|^2_{L^2(\R^{n_1},L^{q}(\R^{n_2}))}$ uniformly in $j$, the above sum is bounded by a constant times
\begin{align}
  \sum_{i'}2^{2s_1i'}\|\varphi_{i'}g\|^2_{L^2(\R^{n_1},L^{q}(\R^{n_2}))}\sum_{j < i' +3}2^{2s_2j}\|\Phi_{i'-3}\psi_jf\|_{L^\infty(\R^{n_1},L^p(\R^{n_2}))}^2. \label{sumpq}
\end{align}
For sufficiently large $t_2$ and $t_2'$ we have embeddings
\begin{align}
  H^{t_2,2}(\R^{n_2}) &\subseteq L^p(\R^{n_2})\label{t2}\\
  H^{t_2',2}(\R^{n_2}) & \subseteq L^q(\R^{n_2}).\label{t2'}
\end{align}
Thus, we can bound the sum (\ref{sumpq}) by
\begin{align}
  \sum_{i'}2^{2s_1i'}\|\varphi_{i'}g\|_{L^2(\R^{n_1},H^{t_2',2}(\R^{n_2}))}^2\sum_{j < i'+3}2^{2s_2j}\|\Phi_{i'-3}\psi_jf\|_{L^\infty(\R^{n_1},H^{t_2,2}(\R^{n_2}))}^2. \label{piece1}
\end{align}
On the other hand, we also have
\begin{align}
  2^{2s_2j}\|\Phi_{i'-3}\psi_jf\|_{L^\infty(\R^{n_1},H^{t_2,2}(\R^{n_2}))}^2 & \sim 2^{2(s_2+t_2)j}\|\Phi_{i'-3}\psi_jf\|_{L^\infty(\R^{n_1},L^2(\R^{n_2}))}^2 \nonumber\\
  & \leq  2^{2(s_2+t_2)j}\|\Phi_{i'-3}\psi_jf\|_{L^2(\R^{n_2},L^\infty(\R^{n_1}))}^2 \nonumber\\
    & \leq  C\cdot2^{2(s_2+t_2)j}\|\Phi_{i'-3}\psi_jf\|_{L^2(\R^{n_2},H^{\frac{n_1}{2}+\eps,2}(\R^{n_1}))}^2. \label{piece2}
\end{align}
The second line uses Minkowski's inequality, and the other inequalities use Sobolev embedding and the support properties of $\Phi_{i'-3}\psi_jf$.  Here, $\eps > 0$ is arbitrary.  Thus, from (\ref{piece1}) and (\ref{piece2}), we obtain the following bound for (\ref{sum1d}):
\begin{align}
  C\sum_{i'}\|\varphi_{i'}g\|_{SH^{s_1,t_2',2}}^2\sum_{j<i'+3}\|\Phi_{i'-3}\psi_jf\|_{SH^{n_1/2+\eps,s_2+t_2,2}}^2. \label{piece3}
\end{align}
For every $i'$, the sum over $j$ is bounded by $\|f\|_{SH^{n_1/2+\eps,s_2+t_2,2}}^2$, uniformly in $i'$.  Thus, we can bound (\ref{piece3}) by a constant times
\begin{equation}
  \|g\|_{SH^{s_1,t_2',2}}^2\|f\|_{SH^{n_1/2+\eps,s_2+t_2,2}}^2.
\end{equation}
So that this estimate can be used in the proposition we are trying to prove, then since $f \in H^{(s_1,s_2'),2}$ and $g \in H^{(s_1,s_2''),2}$, we need to have
\begin{align*}
  SH^{n_1/2+\eps,s_2+t_2,2} &\subseteq H^{(s_1,s_2'),2}\\
  SH^{s_1,t_2',2} &\subseteq H^{(s_1,s_2''),2}.
\end{align*}
This means we must have
\begin{align*}
  s_1 & > n_1/2\\
  s_2+t_2 &< s_1+s_2'-n_1/2\\
  t_2' &\leq s_2'',
\end{align*}
subject to the requirements (\ref{t2}), (\ref{t2'}), and $1/p+1/q=2$. Likewise, when we go to Case 3, the roles of $s_2'$ and $s_2''$ become reversed, and so we also need to satisfy
\begin{align*}
  s_2+t_2 &< s_1+s_2''-n_1/2\\
  t_2' &\leq s_2'.
\end{align*}
Simple arithmetic shows that $s_2 > 0$ satisfying the hypothesis of the theorem ensures that all the above constraints are met.
\end{proof}

\begin{Remark}
  Let $f = \varphi_i(\xi^{(1)})\psi_j(\xi^{(2)})$ and $g = \varphi_{i'}(\xi^{(1)})\psi_{j'}(\xi^{(2)})$ be product bump functions supported on $[2^{i-2},2^{i+2}] \times [2^{j-2},2^{j+2}]$ and $[2^{i'-2},2^{i'+2}]\times[2^{j'-2},2^{j'+2}]$ in frequency space respectively.  We assume the $i,j,i',j'$ satisfy $j' \lesssim j \lesssim i'$ and $i \lesssim i'$ as in Case 1(d) in the above proposition.  Then an elementary calculation shows that
  \begin{align*}
    \|f\|_{H^{(s_1,s_2'),2}}^2 & \sim 2^{2s_1\max(i,j)+2s_2'j}\cdot 2^{in_1}\cdot 2^{jn_2}\\
    \|g\|_{H^{(s_1,s_2''),2}}^2 & \sim 2^{2s_1i' + 2s_2''j'}\cdot 2^{i'n_1}\cdot 2^{j'n_2}\\
    \|fg\|_{H^{(s_1,s_2''),2}}^2 & \sim 2^{2s_1i'+2s_2j}\cdot 2^{(2i + i')n_1} \cdot 2^{(2j'+j)n_2}.
  \end{align*}
  The only way for the last expression to be bounded uniformly in terms of the product of the first two expressions is for $s_2 \leq s_1+s_2'+s_2'' - (n_1+n_2)/2$. This suggests that the conclusion of the above proposition may be close to optimal, though it is unclear to the author how to turn this into a precise sharpness result since the intersection of an anisotropic space with $L^\infty$ is a rather mysterious space.
\end{Remark}

\begin{Corollary}
  If $s_1 > n_1/2$ and $s_1+s_2 > (n_1+n_2)/2$, then $H^{(s_1,s_2),2}(\R^{n_1}\times\R^{n_2})$ is a Banach algebra.
\end{Corollary}

\Proof We have $H^{(s_1,s_2),2}(\R^{n_1}\times\R^{n_2}) \hookrightarrow C^0(\R^{n_1}, H^{s_1-n_1/2+s_2}(\R^{n_2})) \hookrightarrow C^0(\R^{n_1}, C^0(\R^{n_2}))$ by Sobolev embedding.  Hence
$H^{(s_1,s_2),2}(\R^{n_1}\times\R^{n_2}) \subseteq L^\infty(\R^{n_1}\times\R^{n_2})$ and the previous proposition implies the result.\End

\section{Elliptic Boundary Value Problems}

Having established the basic properties of the spaces $A^{(s_1,s_2),p}$ in the previous sections, we want to apply them to the study of elliptic boundary value problems.  Namely, we want to consider elliptic operators acting between sections of vector bundles over a compact manifold with boundary, where we topologize the space of sections with the anisotropic function space topologies.  To do this, we first have to define our anisotropic function spaces on (products of) compact manifolds.  However, this is automatic from the definition of the anisotropic function spaces on open subsets of Euclidean space via the use of local charts (where the charts on a product of manifolds are obtained via the product of the individual charts).  Lemma \ref{LemmaBasic} tells us that different choices of charts yield equivalent norms.

Thus, let $X$ be a compact manifold and let $E$ be a vector bundle over $X$.  Henceforth, we assume all vector bundles are equipped with an inner product, so that in the usual way, function space topologies on $X$ induce corresponding ones for the space of sections of vector bundles over $X$.  Let $C^\infty(E)$ denote the space of smooth sections of $E$ and let $A^{s,p}(E)$ denote the closure of this space with respect to the $A^{s,p}(X)$ topology.  When $X = X_1 \times X_2$, we likewise define $A^{(s_1,s_2),p}(E)$ to be the closure of $C^\infty(E)$ in the $A^{(s_1,s_2),p}(X_1\times X_2)$ topology.  Suppose $\partial X$ is nonempty and let $E_{\partial X} = E|_{\partial X}$ denote the restriction of $E$ to the boundary.  If we are given an elliptic differential operator $A: C^\infty(E) \to C^\infty(F)$ of order $m$ mapping the space of smooth sections of the vector bundle $E$ to those of the vector bundle $F$, we can consider boundary conditions for the operator $A$ as follows.  We have the $(m-1)$th order trace map
\begin{align}
r_{m-1}: C^\infty(E) \to \oplus_{j=0}^{m-1} C^\infty(E_{\partial X})\\
u \mapsto (u, \partial_t u, \ldots, \partial_t^{m-1}u)|_{\partial X},
\end{align}
where $t$ is in the inward normal coordinate in a collar neighborhood of $\partial X \subset X$.   A (pseudodifferential) boundary condition $B$ is a pseudodifferential map from $\oplus_{j=0}^{m-1}C^\infty(E_{\partial X}) \to C^\infty(V)$, where $V$ is some vector bundle over $\partial X$.  Given such a boundary condition, we get two associated operators
\begin{align}
(A,Br_{m-1}): C^\infty(E) &\to C^\infty(F) \oplus C^\infty(V) \label{(A,B)} \\
A_B: \{u \in C^\infty(E): Br_{m-1}u = 0\} &\to C^\infty(F). \label{A_B}
\end{align}
Thus, the first operator is the full mapping pair associated to $A$ and $B$ and the second operator $A_B$ is the operator $A$ restricted to the space of configurations which satisfy the boundary condition $Br_{m-1}u=0$.

The fundamental problem concerning the pair of an elliptic differential operator $A$ and a boundary condition $B$ is to determine what restrictions on $B$ ensure that the operators (\ref{(A,B)}) and (\ref{A_B}) yield Fredholm operators (after taking suitable function space completions). Stated another way, we want to know when the operators (\ref{(A,B)}) and (\ref{A_B}) satisfy a priori elliptic estimates.

When $B$ is a differential operator (i.e., $B$ is a local boundary condition), the study of elliptic boundary value problems as described above go back to the classical work of \cite{ADN}.  There, the function spaces considered are the classical Sobolev spaces $H^{k,p}$ where $k$ is a nonnegative integer and $1<p<\infty$.  Thus, for example, for $A$ and $B$ satisfying suitable hypotheses, one obtains from \cite{ADN} the following elliptic estimate
\begin{align}
\|u\|_{H^{k+m,p}} \leq C(\|Au\|_{H^{k,p}}  + \|u\|_{H^{k,p}}) \label{ellest}
\end{align}
for all $u$ such that $Br_{m-1}u = 0$.  One also has analogous estimates for the full mapping pair $(A,B)$.  Subsequent results were also obtained by H\"ormander \cite{Ho} and Seeley \cite{Se2} for when $B$ is a pseudodifferential operator satisfying certain conditions on its symbol.  Moreover, because pseudodifferential operators are bounded on a variety of function spaces, in particular, Bessel potential and Besov spaces, the methods of H\"ormander and Seeley apply to the scale of spaces $H^{s,p}$ and $B^{s,p}$.  One can also work with the more general scale of spaces $F^s_{p,q}$ and $B^s_{p,q}$, see \cite{RS}.

For us, we will be interested in the case when $X = X_1 \times X_2$ is a product and we wish to obtain elliptic estimates on the anisotropic spaces.  Here, we assume $X_1$ is a compact manifold with boundary and $X_2$ is a closed manifold, so that the anisotropic function spaces on $X$ and $\partial X$ are defined with respect to the decompositions $X = X_1\times X_2$ and $\partial X = \partial X_1 \times X_2$.  If we follow the approach of \cite{Se2}\footnote{Other approaches are possible, but to the author's knowledge, Seeley's is the cleanest and most concise.}, then it is a simple matter to verify that Seeley's methods generalize to the anisotropic function spaces $A^{(s_1,s_2),p}$.  This is because Seeley constructs a parametrix for the boundary value problems considered, and these parametrices, being compositions of pseudodifferential, trace, and extension operators, map anisotropic spaces to each other from the results of the previous sections.  Here, it is key that the anisotropy is tangential to the boundary; this way, the anisotropy is preserved under taking traces and extensions when passing between $X$ and $\partial X$.  From this, one can reprove much of the results of \cite{Se2} in the context of anisotropic spaces.  When $p=2$, the spaces $H^{(s_1,s_2),2}$ and their application in elliptic boundary value problems were already explicitly considered in H\"ormander \cite{Ho} (see also \cite{WRL}). Our results are the direct generalization to $p \neq 2$.

In the remainder of this section, we wish to check the above statements more precisely under slightly more general notions of a boundary condition.  As mentioned above, when $B$ is a pseudodifferential operator, there is a standard condition on its principal symbol which defines $B$ to be an \textit{elliptic boundary condition} (see \cite{Ho}).  In \cite{Se2}, Seeley defines a slightly more general notion of a \textit{well-posed} boundary condition (also called \textit{injectively elliptic} in the literature), for which he establishes elliptic estimates for the associated boundary value problems.  However, the most general notion of a boundary condition appearing the literature, to the author's knowledge, is stated in \cite{SSS} (although such a notion is already obviously implicit in the work of \cite{Se2}).  Namely, consider the following situation, which is the situation in \cite{SSS} restated for anisotropic function spaces.  For $s_1,s_2 \geq 0$, we have the operator\footnote{For specificity, we work with anisotropic Bessel potential spaces $H^{(s_1,s_2),p}$ on $X$, though everything we do applies equally well to the anisotropic Besov spaces $B^{(s_1,s_2),p}$. Moreover, the methods of \cite{Se2} also allow for $s_1 < 0$ under natural restrictions on the order of $B$.  To keep matters simple, we will consider $s_1 \geq 0$ since the corresponding results for $s_1 < 0$ easily follow through.}
\begin{equation}
  A: H^{(s_1+m,s_2),p}(E) \to H^{(s_1,s_2),p}(F). \label{Aop}
\end{equation}
The $(m-1)$th order trace map gives us a bounded operator
\begin{equation}
 r_{m-1}: H^{(s_1+m,s_2),p}(E) \to \oplus_{j=0}^{m-1} B^{(s_1+m-1/p-j,s_2),p}(E_{\partial X}).
\end{equation}
Here, we apply Theorem \ref{ThmATrace} in the manifold setting.  Of course, all the results of the previous sections on Euclidean space easily generalize in the appropriate way to manifolds, since manifolds are locally Euclidean.

\begin{Definition}\label{DefBC}
  A boundary condition $B$ for the $m$th order elliptic operator $A$ given by (\ref{Aop}) is a map
\begin{equation}
  B: \oplus_{j=0}^{m-1} B^{(s_1+m-1/p-j,s_2),p}(E_{\partial X})\to \X \label{BC}
\end{equation}
from the boundary data space to some Banach space $\X$.
\end{Definition}

In the typical situation where $B$ is a (pseudo)differential operator acting between sections of vector bundles over $\partial X$, then $\X$ is explicitly determined by $B$ acting on the various components of $\oplus_{j=0}^{m-1} B^{(s_1+m-1/p-j,s_2),p}(E_{\partial X})$.  For instance, if $A$ is the Laplacian and $B$ is the Dirichlet or Neumann boundary condition, then $B$ corresponds to the projection onto the first or second factor, respectively, of the image of $r_1$.  However, the point of the Definition \ref{DefBC} is that $B$ may be any abstract map of Banach spaces, and as well shall now see, its merely the formal properties of $B$ that are needed in order to obtain the usual elliptic estimates one gets when $B$ is a pseudodifferential elliptic boundary condition.

For any $s_1 \in \R$ and $s_2 \geq 0$, consider the space
$$Z^{(s_1,s_2),p}(A) := \ker A \subset H^{(s_1,s_2),p}(E),$$
the elements of $H^{(s_1,s_2),p}(E)$ which lie in the kernel of $A$.  For $s_2 = 0$, it is proven in \cite{Se1} that we have a bounded map
\begin{equation}
  r_{m-1}: Z^{(s_1,s_2),p}(A) \to \oplus_{j=0}^{m-1} B^{(s_1-1/p-j,s_2),p}(E_{\partial X}) \label{rZ}
\end{equation}
for all $s_1 \in \R$.  This implies (\ref{rZ}) is bounded for all $s_2 > 0$ by Theorem \ref{ThmATrace}.  In addition, when $s_2 = 0$, the results of \cite{Se1,Se2} tell us that there exist special pseudodifferential operators which intertwine $Z^{s_1,p}(A)$ with its Cauchy data on the boundary.  Namely, we have a pseudodifferential projection $P^+$ which maps the boundary data space $\oplus_{j=0}^{m-1} B^{s_1-1/p-j,p}(E_{\partial X})$ onto $r_{m-1}(Z^{s_1,p}(A))$.  Furthermore, we have a map $P$ from the boundary data space $\oplus_{j=0}^{m-1} B^{s_1-1/p-j,p}(E_{\partial X})$ to the nullspace $Z^{s_1,p}(A)$ such that $r_{m-1}P = P^+$ and $PP^+ = P$.  Furthermore, $\im P \subset Z^{s_1,p}(A)$ is complemented by the space $Z_0(A) = \{u \in C^\infty(E) : Au = 0, r_{m-1}u=0\}$, which is finite dimensional.  Moreover, smooth configurations are dense in $\im P^+$ and $Z^{s_1,p}(A)$.  The operators $P^+$ and $P$ are known as a \textit{Calderon projection} and \textit{Poisson operator}, respectively.  Note that these operators are not unique, since although their ranges are specified, their kernels are not.

The following is an anisotropic generalization of these results.

\renewcommand{\theenumi}{(\roman{enumi})}
\begin{Lemma}\cite{Se1}\label{LemmaSeeley} Let $s_1 \in \R$, $s_2 \geq 0$, $1<p<\infty.$
\begin{enumerate}
\item The Calderon projection $P^+: \oplus_{j=0}^{m-1} B^{(s_1-1/p-j,s_2),p}(E_{\partial X}) \to r_{m-1}(Z^{(s_1,s_2),p}(A))$ is bounded.
\item The Poisson operator $P^+: \oplus_{j=0}^{m-1} B^{(s_1-1/p-j,s_2),p}(E_{\partial X}) \to Z^{(s_1,s_2),p}(A)$ is bounded.
\end{enumerate}
\end{Lemma}

\begin{proof}
 (i) The Calderon projection is given by a (matrix of) pseudodifferential operators.  By Theorem \ref{ThmProdPSDO}, since pseudodifferential operators are bounded on anisotropic function spaces, the result follows.

 (ii) In the isotropic case $s_2 = 0$, the statement follows from \cite[Lemma 4]{Se1}.  There, the statement is proven in three steps: first for $s_1 < 1/p$, second for $s_1$ an integer greater than or equal to $m$, and then finally, the result follows for all $s_1$ by interpolation.  By the interpolation result, Lemma \ref{LemmaInterp}, to prove (ii), it suffices to check that Seeley's first two steps also hold for any fixed $s_2 \geq 0$.  However, because traces, extensions, and pseudodifferential operators preserve tangential anisotropy by our previous results, one can check straightforwardly that Seeley's proofs carry through unchanged.
 \end{proof}

We now have our main theorem for the full mapping pair, which generalizes \cite[Theorem 2]{SSS} to anisotropic function spaces.

\begin{Theorem}\label{Thm(A,B)}
 Let $s_1,s_2 \geq 0$, $1<p<\infty$, and let $A$ be an $m$th order elliptic differential operator.  Let $B$ be a boundary condition for $A$ as in (\ref{BC}). Then the operator $(A,Br_{m-1}): H^{(s_1+m,s_2),p}(E) \to H^{(s_1,s_2),p}(F) \oplus \X$ is Fredholm if and only if  $B: \im P^+ \to \X$ is Fredholm.
\end{Theorem}

\begin{proof}
 The proof proceeds as in \cite[Theorem 2]{SSS}, which closely follows the methods of \cite{Se2}, only we have to check that the steps follow through for anisotropic function spaces.  For this, we make use of Theorems \ref{ThmProdPSDO} and \ref{ThmATrace} and Lemma \ref{LemmaSeeley}.

 In detail, given $f \in H^{s_1,p}(F)$ which is orthogonal to the finite dimensional space $Z_0(A^*)$, where $A^*$ is the formal adjoint of $A$, there exists  $v \in A^{s_1+m,p}(E)$ such that $Av = f$.  This $v$ is obtained via the following.  First, we construct an ``invertible double" $D$ for $A$ as in \cite{Se2}. Namely, we embed $A$ as a summand of an invertible $m$th order elliptic operator $D$ acting on sections of $\tilde E \oplus \tilde F$, where $\tilde E$ and $\tilde F$ are bundles on a closed manifold $\tilde X \supset X$ which extend the vector bundles $E$ and $F$ over $X$.  We then apply the inverse $D^{-1}$ of the invertible double to $(0, E_k f)$, where $E_k$ is an extension map from $X$ to $\tilde X$ as in (\ref{Ek}) such that $k > s_1$ and all functions in the image of $E_k$ have support in some small tubular neighborhood of $X \subset \tilde X$.  Our configuration $v$ is then the restriction of $D^{-1}(0, E_k f)$ to the bundle $E$ over $X$.  Since the extension map $E_k$ and pseudodifferential operators preserve tangential anisotropy, it follows that if $f \in H^{(s_1,s_2),p}(F)$, we also have $v$ as constructed above lies in $H^{(s_1+m,s_2),p}(E)$.

To establish the theorem, we consider the generalized boundary value problem
\begin{align}
  Au &= f \in A^{(s_1,s_2),p}(F) \label{BVPeq1}\\
  Br_{m-1}u &= g \in \X. \label{BVPeq2}
\end{align}
We can assume $f$ is orthogonal to $Z_0(A^*)$, since this is only a finite dimensional restriction, which does not affect the Fredholm properties of the full mapping pair $(A,B)$.  Then from the previous steps, there exists a $v \in H^{(s_1+m,s_2),p}(E)$ such that $Av = f$.  Thus, $A(u-v) = 0$, and letting $w = u-v$, we are now reduced to solving the problem
\begin{align}
  Br_{m-1}w &= g - Br_{m-1}v \in \X, \qquad w \in \im P^+. \label{BVPeq3}
\end{align}
Since every $f \in H^{(s_1,s_2),p}(F)$ in the complement of the finite dimensional space $Z_0(A^*)$ determines a unique $v \in H^{(s_1+m,s_2),p}(E)$ as described above, it is clear that the problem (\ref{BVPeq1})-(\ref{BVPeq2}) is Fredholm if and only if $B: \im P^+ \to \X$ is Fredholm.
\end{proof}

Typically, the boundary condition $B$ is given by a pseudodifferential operator with closed range, in which case, the standard definition of an elliptic boundary condition (in terms of the principal symbol of $B$) ensures that $B: \im P^+ \to \im B$ is Fredholm.  When $B$ is differential and its range is equal to a Banach space of sections of a vector bundle, such a condition is also known as the Lopatinsky-Shapiro condition.

From Theorem \ref{Thm(A,B)}, we get a corresponding theorem on the restricted operator $A_B$ associated to the full mapping pair $(A,B)$.  Namely, consider the restricted domain
$$H^{(s_1+m,s_2),p}_B(E) = \{u \in H^{(s_1+m,s_2),p}(E) : Br_{m-1}u = 0\}$$
and the operator
\begin{equation}
  A_B: H^{(s_1+m,s_2),p}_B(E) \to H^{(s_1,s_2),p}(F). \label{opA_B}
\end{equation}

\begin{Theorem}\label{ThmA_B}
Let $s_1,s_2\geq 0$ and $1<p<\infty$. We have the following:
\begin{enumerate}
\item The operator (\ref{opA_B}) is Fredholm if and only if $B: \im P^+ \to \X$ is Fredholm.
\item The kernel of (\ref{opA_B}) is spanned by finitely many smooth configurations if and only if $\ker B \cap \im P^+$ is spanned by finitely many smooth configurations.
\item The range of (\ref{opA_B}) is complemented by the span of finitely many smooth configurations if and only if $\ker B + \im P^+$ is complemented by the span of finitely many smooth configurations.
\end{enumerate}
\end{Theorem}

\begin{proof}
(i) Suppose the map $B: \im P^+ \to \X$ is Fredholm.  This means that the quotient of the boundary value space  $\oplus_{j=0}^{m-1} B^{(s_1+m-1/p-j,s_2),p}(E_{\partial X})$ by the closed subspace $\ker B$ is, modulo finite dimensional subspaces, isomorphic to $\im P^+$.  (When we write $\ker B$, we always mean the kernel of $B$ acting on all of $\oplus_{j=0}^{m-1} B^{(s_1+m-1/p-j,s_2),p}(E_{\partial X})$.) More precisely, $B: \im P^+ \to \X$ is Fredholm if and only if $\ker B \cap \im P^+$ is finite dimensional and $\ker B + \im P^+$ is closed and complemented by a finite dimensional subspace.  Given $f \in H^{(s_1,s_2),p}(F)$ orthogonal to $Z_0(A^*)$, we can find a $v \in H^{(s_1+m,s_2),p}(E)$ as in Theorem \ref{Thm(A,B)} such that $Av = f$.  By the above, for $f$ lying in the complement of some finite dimensional subspace of $H^{(s_1,s_2),p}(F)$, such a $v$ will satisfy $r_{m-1}v \in \ker B + \im P^+$.  Moreover, there exists a unique projection $\pi: \ker B + \im P^+ \to \ker B$ with $\ker \pi \subseteq \im P^+$.  It follows that we can write $r_{m-1}v = x + x'$, where $x = \pi(r_{m-1}v)$ and $x' = (1 - \pi)(r_{m-1}v)$.  Let $v' := v - Px'$, where $P$ is the Poisson operator of $A$ as in Lemma \ref{LemmaSeeley}.  Then $A(v - Px') = Av = f$, since $\im P \subset \ker A$.  Moreover, $r_{m-1}(v - Px') = r_{m-1}v - x' \in \ker B.$  In this way, we have constructed a $v' \in H^{(s_1+m,s_2),p}_B(E)$ such that $Av' = f$, whenever $f$ lies in some subspace of finite codimension in $H^{(s_1,s_2),p}(F)$.  Moreover, if $Av'' = f$ with $v'' \in H^{(s_1+m,s_2),p}_B(E)$, then $v' - v'' \in \ker A$ and so $r(v' - v'')$ belongs to the finite dimensional space $\ker B \cap \im P^+$.  It follows that $\ker A_B$ is isomorphic to the finite dimensional space $(\ker B \cap \im P^+) \oplus Z_0(A)$.  Thus, the kernel and cokernel of (\ref{opA_B}) are finite dimensional, which shows that (\ref{opA_B}) is Fredholm.

Similar reasoning leads to the reverse conclusion.  Namely, suppose (\ref{opA_B})is Fredholm.  Then since $A_B$ has finite dimensional kernel, we have $\ker B \cap \im P^+$ is finite dimensional.  Moreover, since $A_B$ has closed range and finite dimensional cokernel, repeating the arguments in the previous paragraph shows that $\ker B + \im P^+$ is closed and of finite codimension.  This implies $B: \im P^+ \to \X$ is Fredholm.

(ii),(iii) These statements now easily follow from the proof of (i).  Indeed, one relates the smoothness of elements of the kernel $A_B$ with the smoothness of their boundary values via the Poisson operator $P$.  Likewise, one does a similar thing for elements of the cokernel of $A_B$, whereby one uses the invertible double $D$ to relate smooth elements in the complement of $Z_0(A^*)$ to smooth elements of $H^{(s_1+m,s_2),p}(E)$ and hence to smooth elements of $\oplus_{j=0}^{m-1} B^{(s_1+m-1/p-j,s_2),p}(E_{\partial X})$ via the restriction map $r_{m-1}$.
\end{proof}

Observe that when (i) and (ii) hold above, we obtain the elliptic estimate (\ref{ellest}) for $A_B$, since the last, lowest order term in (\ref{ellest}) is only needed to control the kernel of $A_B$, which from (ii), is a finite dimensional space spanned by smooth sections.  Moreover, the above theorem shows that the elliptic estimate (\ref{ellest}) generalizes to anisotropic function spaces, and it furthermore, it yields an elliptic estimate for the full mapping pair $(A,B)$.  Summarizing, we have

\begin{Corollary}
  Suppose the hypotheses of (i) and (ii) of Theorem \ref{ThmA_B} hold.  Then we have the elliptic estimate
  $$\|u\|_{H^{(s_1+m,s_2),p}(E)} \leq C(\|Au\|_{H^{(s_1,s_2),p}(F)} + \|Br_{m-1}u\|_{\oplus_{j=0}^{m-1} B^{(s_1+m-1/p-j,s_2),p}(E_{\partial X})} + \|\pi u\|),$$
where $\pi$ is any projection onto the kernel of $A_B$ and $\|\cdot\|$ is any norm on the finite dimensional space $\im \pi$.
\end{Corollary}


\noindent\textsc{Massachusetts Institute of Technology, Cambridge, MA 02139}\\
\textit{Email address:} \texttt{timothyn@math.mit.edu} 

\end{document}